\documentclass[12pt]{article}

\usepackage{amsmath,amssymb,amsthm,amsfonts,diagrams}
\usepackage{nicematrix,tikz}
\usetikzlibrary{fit}
\tikzset{highlight/.style={rectangle,
fill=gray!50,
rounded corners = 0.5 mm, 
inner sep=1pt,
fit=#1}}

\def\Hom{{\rm Hom}}
\def\ind{{\rm ind}}
\def\ad{\mathop{\rm ad}}
\def\im{\mathop{\rm im}}

\def\phi{\varphi}
\def\g{\mathfrak g}

\def\h{\mathfrak h}
\def\F{\mathbb F}

\makeatletter \let\@@pmod\pmod
\DeclareRobustCommand{\pmod}{\@ifstar\@pmods\@@pmod}
\def\@pmods#1{\mkern4mu({\operator@font mod}\mkern 6mu#1)}
\makeatother

\makeatletter
\def\foo#1\endgraf\unskip#2\foo{\def\row@to@buffer{#1\endgraf\unskip\unskip#2}}
\expandafter\foo\row@to@buffer\foo
\makeatother

\newtheorem{theorem}{Theorem}[section]
\newtheorem{lemma}[theorem]{Lemma}
\newtheorem{proposition}[theorem]{Proposition}
\newtheorem{corollary}[theorem]{Corollary}
\newtheorem{definition}[theorem]{Definition}
\newtheorem{rem}[theorem]{Remark}

\theoremstyle{remark}

\providecommand{\keywords}[1]{\noindent{Keywords:} #1}
\providecommand{\classify}[1]{\noindent{Mathematics Subject
    Classification:} #1}

\title{Restricted Cohomology of Heisenberg Lie Algebras}

\author{
Yong Yang\\
College of Mathematics and System Science\\
Xinjiang University\\
Urumqi 830046, China\\
yangyong195888221@163.com
}

\date{}

\begin{document}
\maketitle

\begin{abstract}
The Heisenberg Lie
  algebras over an algebraically closed field $\F$ of
  characteristic $p>0$ always admit a family of restricted structures. We use the ordinary 1- and 2-cohomology spaces with
 adjoint coefficients to compute the restricted 1- and 2-cohomology
  spaces of these restricted Heisenberg Lie algebras.  
\end{abstract}

{\footnotesize

\classify{17B50; 17B56}}

\section{Introduction}
Given a symplectic vector space,  its Heisenberg Lie algebra arising from the
symplectic form  is a two-step nilpotent Lie algebra with one-dimensional
center, as shown in the following definition:
\begin{definition}\rm
  For a non-negative integer $m$, the \emph{Heisenberg Lie algebra}
  $\mathfrak{h}_{m}$ is the $(2m+1)$-dimensional vector space over a
  field $\F$ spanned by the elements
  \[\{e_1,\ldots , e_{2m},e_{2m+1}\}\] with the non-vanishing Lie
  brackets
  \[[e_i,e_{m+i}]=e_{2m+1},\] for $1\le i\le m$.
\end{definition}
Due to the  applications in the
commutation relations in quantum mechanics, there are a variety of algebraic works on Heisenberg Lie algebras.
For instance the
cohomology with trivial coefficients of Heisenberg Lie algebras over
fields of characteristic zero was computed  \cite{S},
the algebraic
Morse theory was applied to compute the Poincar\'e polynomial of the Lie algebra
$\mathfrak{h}_m$ over fields of characteristic $p=2$ \cite{SK} and
the $n$-th Betti number of
$\mathfrak{h}_m$ over a field of any prime characteristic for $n \le m$ was given in 
\cite{CJ}.
There are also other results related to invariants of Heisenberg algebras, like gradings and symmetries  \cite{T} and the dimensions of the adjoint cohomology spaces over fields of characteristic zero \cite{A,CaTi,M}.
% adjoint cohomology, see \cite{M, CaTi, CaTi2, A}.

The notion of restricted Lie algebras over fields of positive
characteristic was introduced by Jacobson in 1937 \cite{J}. Since
then, the study of restricted Lie algebras has proved to be fruitful
for several reasons, see \cite{Se} for example.
By Jacobson's Theorem,
 two-step nilpotent Lie algebras are restrictable.
 As a result, the Lie algebra $\h_m$ always has a family of restricted structures parameterized by $\lambda=(\lambda_1,\dots, \lambda_{2m+1})$.  Denote the
corresponding restricted Lie algebra by $\mathfrak{h}_m^\lambda$.  The trivial restricted cohomology and central extensions of $\mathfrak{h}_m^\lambda$ were studied in \cite{EFY}.

Much less is known about the adjoint restricted cohomology of $\h^\lambda_m$.
Our goal in this note is to  compute the 
restricted cohomology spaces 
$H_*^q(\mathfrak{h}_m^{\lambda},\mathfrak{h}_m^{\lambda})$ for $q=1,2$, and give explicit
bases for these spaces. 

\section{Preliminaries and Main Results}

In this section, we recall the restricted Heisenberg Lie algebra, 
the Chevalley-Eilenberg cochain complex,
and the restricted cochain complex \cite{EFu}. Everywhere in this
section, $\F$ denotes an algebraically closed field of characteristic
$p>0$.  For $j\geq 2$ and $g_1,\ldots,g_ j\in \h_m$, we denote the $j$-fold bracket
$[g_1,g_2,g_3,\ldots,g_j]=[[\ldots[[g_1,g_2],g_3],\ldots],g_j].$

\subsection{The Restricted Lie Algebra  $\mathfrak{h}^{\lambda}_m(p)$}
Since $\mathfrak{h}_m$ is 2-step nilpotent, $(\ad g)^p=0$ for all
$g \in \mathfrak{h}_m$ so we can use Jacobson's theorem \cite{J} and
define a restricted Lie algebra $[p]$-map on $\mathfrak{h}_m$ by
choosing $e_i^{[p]}\in \mathfrak{h}_m$ so that
$\ad e_i^{[p]}=(\ad e_i)^p=0$ for all $1\le i\le 2m+1$. That is, by
choosing $e_i^{[p]}$ in the center $\F e_{2m+1}$. For each $i$, we
choose a scalar $\lambda_i\in \F$ and set
$e_i^{[p]}=\lambda_i e_{2m+1}$. Let
$\lambda=(\lambda_1,\dots, \lambda_{2m+1})$, and denote the
corresponding restricted Lie algebra by
$\mathfrak{h}_m^\lambda(p)$.  Moreover, 
For $g=\sum^{2m+1}_{i=1} a_i e_i\in \mathfrak{h}_m^{\lambda}$, we have \cite{EFY}:
\begin{equation}\label{r}
  g^{[p]}  = \left\{
                             \begin{array}{ll}
                            \left(\sum^{2m+1}_{i=1}a^{p}_i\lambda_i\right) e_{2m+1}, & \hbox{$p>2$;}\\
                             \left(\sum^{2m+1}_{i=1}a^{2}_i\lambda_i+\sum^{m}_{j=1} a_j
     a_{m+j} \right)e_{2m+1}. & \hbox{$p=2$.}
                              
                             \end{array}
                                                                                   \right.
\end{equation}
 Below we write $\h_m^\lambda$ in place of
$\h_m^\lambda(p)$ when no confusion can arise.
For $i: 1\leq i\leq 2m+1$, we define the  conjugate $i'$  as follows:
$$
  i’ = \left\{
                             \begin{array}{lll}
                           i+m, & \hbox{$1\leq i\leq m$;}\\
                            i-m, & \hbox{$m+1\leq i\leq 2m$;}\\
                               2m+1, & \hbox{$i=2m+1$.}
                             \end{array}
                                                                                   \right.
$$

%From now on, $\mathfrak{g}$ denotes a restricted Lie algebra.

\subsection{Cochain Complexes with Adjoint Coefficients}

\subsubsection{Ordinary Cochain Complex}
We here  only describe the
Chevalley-Eilenberg cochain spaces $C^q(\h_m,\h_m)$
for $q=0,1,2,3$ and differentials $d^q:C^q(\h_m,\h_m)\to C^{q+1}(\h_m,\h_m)$
for $q=0,1,2$. For details on the Chevalley-Eilenberg cochain
complex we refer the reader to \cite{ChE,F}. Set
$C^0 (\h_m,\h_m)=\h_m$ and $C^q (\h_m,\h_m)= \mathrm{Hom}(\wedge^q\h_m,\h_m)$ for $q=1,2,3$. We will use the
following bases throughout the paper.
\begin{align*}
  C^0 (\h_m,\h_m):&\ \{e_i\ |\ 1\le i\le 2m+1\}&\\
  C^1 (\h_m,\h_m):&\  \{e_j^i\ |\ 1\le i,j\le 2m+1\}&\\
  C^2 (\h_m,\h_m):&\ \{e_k^{i,j}\ |\ 1\le i<j\le 2m+1,1\le k\le 2m+1\}&\\
  C^3 (\h_m,\h_m):&\ \{e_{l}^{i,j,k}\ |\ 1\le i<j<k\le 2m+1,1\le l \le 2m+1\}&\\
\end{align*}
Here $e^i_j: e_i\rightarrow e_j$, $e^{i,j}_k: (e_i,e_j)\rightarrow e_k$ and $e_{l}^{i,j,k}: (e_i,e_j,e_k)\rightarrow e_{l}$ denote  the
basis vectors, respectively.
The differentials $d^q:C^q (\h_m,\h_m)\to C^{q+1}(\h_m,\h_m)$ are defined for $e_q,e_r,e_s\in\h_m$ by\small
\begin{align*}
  &d^0: C^0 (\h_m,\h_m)\to C^1 (\h_m,\h_m),  \  d^0(e_i)(e_q)=[e_i,e_q]&\\
 & d^1:C^1 (\h_m,\h_m)\to C^2 (\h_m,\h_m), \   d^1(e^i_j)(e_q\wedge e_r)=[e_q,e_j^i(e_r)]-[e_r,e^i_j(e_q)]-e^i_j([e_q,e_r])&\\
 & d^2:C^2 (\h_m,\h_m)\to C^3 (\h_m,\h_m), \  d^2(e^{i,j}_k)(e_q\wedge e_r\wedge e_s)=[e_q,e^{i,j}_k(e_r,e_s)]-[e_r,e^{i,j}_k(e_q,e_s)]
\\
&+[e_s,e^{i,j}_k(e_q,e_r)]-e^{i,j}_k([e_q,e_r]\wedge e_s)+e^{i,j}_k([e_q,e_s]\wedge e_r)-e^{i,j}_k([e_r,e_s]\wedge e_q). &\\
\end{align*}\normalsize
The maps $d^q$ satisfy $d^{q}d^{q-1}=0$ and 
$H^q(\h_m,\h_m)=\ker(d^q)/\im(d^{q-1})$.  
From now on, for a cocycle $\phi\in \ker(d^q)$, we use $\phi$ to represent both the cocycle and the corresponding cohomology class when no confusion can arise.

\subsubsection{Restricted Cochain Complex}

In this subsection, we recall the definitions and results on the (partial)
restricted cochain complex given in \cite{EFu} for $\h_m^{\lambda}$ in the case of adjoint
coefficients. For a more general treatment of the cohomology of
restricted Lie algebras, we refer the reader to \cite{EFu,Fe}.

Given $\phi\in C^2(\h_m,\h_m)$, a map $\omega:\h_m\to\h_m$ is {\bf
  $\phi$-compatible} if for all $g,h\in\h_m$ and all $a\in\F$
\\
  
$\omega(a g)=a^p \omega (g)$ and
\begin{align}
  \label{starprop}
\begin{split}
  \omega(g+h)=&\omega(g)+\omega(h) + \sum_{\substack{g_i=\mbox{\rm\scriptsize $g$
        or $h$}\\ g_1=g, g_2=h}}
  \frac{1}{\#(g)}(\phi([g_1,g_2,\dots,g_{p-1}]\wedge g_p)\\
&-[g_p,\phi([g_1,g_2,\dots,g_{p-2}]\wedge g_{p-1})])\\
\end{split}
\end{align}
where $\#(g)$ is the number of factors $g_i$ equal to $g$. 

We define
\begin{align*}
\Hom_{\rm Fr}(\h_m,\h_m) = \{f:\h_m\to\h_m\ &|\ f(a g+b
  h)=a^pf(g)+b^pf(h)\ \mathrm{for}\  \mathrm{all}\ \\
&a,b\in\F \ \rm{and}\ \mathrm{all}\  \emph{g,h}\in \h_\emph{m}\}
\end{align*}
 to be
the space of {\it Frobenius homomorphisms} on $\h_m$. A map
$\omega:\h_m\to\h_m$ is $0$-compatible  if and only if
$\omega\in \Hom_{\rm Fr}(\h_m,\h_m)$.

Given $\phi\in C^2(\h_m,\h_m)$, similar to the trivial case \cite{EFi2,EFY}, we can get a unique $\phi$-compatible map  $\omega$ by assigning the values of $\omega$ arbitrarily on a
basis of $\h_m$. In
particular, we can define $\tilde\phi(e_i)=0$ for all $i$ and use
(\ref{starprop}) to determine a unique $\phi$-compatible map
$\widetilde\phi:\h_m\to\h_m$. Note that, in general, $\widetilde\phi\ne 0$
but $\widetilde\phi (0)=0$. Moreover, If $\phi_1,\phi_2\in C^2(\h_m,\h_m)$
and $a\in\F$, then
$\widetilde{(a\phi_1+\phi_2)} = a\widetilde\phi_1 + \widetilde\phi_2$.
\\

Given $\zeta\in C^3(\h_m,\h_m)$, a map $\eta:\h_m\times \h_m\to\h_m$ is {\bf
  $\zeta$-compatible} if for all $a\in\F$ and all $g,h,h_1,h_2\in\g$,
$\eta(\cdot,h)$ is linear in the first coordinate,
$\eta(g,a h)=a^p\eta(g,h)$ and
\begin{align*}
  \eta(g,h_1+h_2) &=
                    \eta(g,h_1)+\eta(g,h_2)-\nonumber \\
                  & \sum_{\substack{l_1,\dots,l_p=1 {\rm or} 2\\ l_1=1,
  l_2=2}}\frac{1}{\#\{l_i=1\}}(\zeta (g\wedge
  [h_{l_1},\cdots,h_{l_{p-1}}]\wedge h_{l_{p}})\\
&-[h_{l_p},\zeta (g\wedge
  [h_{l_1},\cdots,h_{l_{p-2}}]\wedge h_{l_{p-1}})]).
\end{align*}
The restricted cochain spaces are defined as 
\begin{align*}
&C^0_*(\h^{\lambda}_m,\h^{\lambda}_m)=C^0 (\h_m,\h_m),\\
&C^1_*(\h^{\lambda}_m,\h^{\lambda}_m)=C^1 (\h_m,\h_m),\\
&C^2_*(\h^{\lambda}_m,\h^{\lambda}_m)=\{(\phi,\omega)\ |\ \phi\in C^2 (\h_m,\h_m), \omega:\h_m\to\h_m\
  \mbox{\rm is $\phi$-compatible}\},\\
&C^3_*(\h^{\lambda}_m,\h^{\lambda}_m)=\{(\zeta,\eta)\ |\ \zeta\in C^3 (\h_m,\h_m),
  \eta:\h_m\times\h_m\to\h_m\ \mbox{\rm is $\zeta$-compatible}\}.
\end{align*}
 % f
% $\phi\in C^2 (\g)$, we can assign values $\omega (e_i)$ arbitrarily
% to basis elements $e_i$ for $\g$. Set
% $\omega(\alpha e_i)=\alpha^p\omega(e_i)$ for all $\alpha\in\F$ and
% use (\ref{starprop}) to determine a unique $\phi$-compatible map
% $\omega:\g\to\F$ \cite{EFi2}.

For $1\le i,j\le 2m+1$, define $\overline e_j^i:\h_m\to\h_m$ by
$\overline e_j^i \left(\sum_{k=1}^{2m+1} a_k e_k\right ) = a_i^pe_j.$ The set
$\{\overline e_j^i\ |\ 1\le i,j\le 2m+1\}$ is a basis for the space of Frobenius homomorphisms
$\Hom_{\rm Fr}(\h_m,\h_m)$. Since the dimension of  the restricted cochain space $C^k(\g,M)$ of a restricted Lie algebra $\g$ with coefficients in $M$  is the same as the tensor product
$S^k\g\otimes M$ \cite{EFu}, we have
\[\dim C^2_*(\h^{\lambda}_m,\h^{\lambda}_m) = \binom{\dim\h_m+1}{2}\dim\h_m=\binom{\dim\h_m}{2}\dim\h_m+(\dim\h_m)^2.\]
Therefore
\begin{equation*}
  \{(e_k^{i,j},\widetilde{e_k^{i,j}})\ |\ 1\le i<j\le 2m+1,1\le k\le 2m+1\} \cup \{(0,\overline e_j^i)\ |\ 1\le i,j\le 2m+1\}
\end{equation*}
is a basis for $C^2_*(\h^{\lambda}_m,\h^{\lambda}_m)$. We use this basis in all computations that
follow.

Define $d_*^0=d^0$. For $\psi\in C^1_*(\h^{\lambda}_m,\h^{\lambda}_m)$, define the map
$\ind^1(\psi):\h^{\lambda}_m\to\h^{\lambda}_m$ by
$
\ind^1(\psi)(g)=\psi(g^{[p]})-(\ad g)^{p-1}(\psi(g)).
$
 The map
$\ind^1(\psi)$ is $d^1(\psi)$-compatible for all $\psi\in C^1_*(\h^{\lambda}_m,\h^{\lambda}_m)$,
and the differential $d^1_*:C^1_*(\h^{\lambda}_m,\h^{\lambda}_m)\to C^2_*(\h^{\lambda}_m,\h^{\lambda}_m)$ is defined by
$ d^1_*(\psi) = (d^1(\psi),\ind^1(\psi)).$

For $(\phi,\omega)\in C^2_*(\h^{\lambda}_m,\h^{\lambda}_m)$, define the map
$\ind^2(\phi,\omega):\h^{\lambda}_m \times\h^{\lambda}_m \to\h^{\lambda}_m$ by the formula
\begin{align}\label{ind2}
\begin{split}
\ind^2(\phi,\omega)(g,h) =&\phi(g\wedge h^{[p]})+[g,\omega(h)]+\phi(h\wedge (\ad h)^{p-1}(g)) \\                       
                         &+[h,\phi(h\wedge(\ad h)^{p-2}(g))].
\end{split}
\end{align}
The map $\ind^2(\phi,\omega)$ is $d^2(\phi)$-compatible for all
$\phi\in C^2(\h_m,\h_m)$, and the differential 
$d^2_*:C^2_*(\h^{\lambda}_m,\h^{\lambda}_m)\to C^3_*(\h^{\lambda}_m,\h^{\lambda}_m)$
is defined by
$
  d^2_*(\phi,\omega) =
  (d^2(\phi),\ind^2(\phi,\omega)). 
$
These maps $d_*^q$ satisfy $d_*^{q}d_*^{q-1}=0$ and we define
\[H_*^q(\h^{\lambda}_m,\h^{\lambda}_m)=\ker(d_*^q)/\im(d_*^{q-1}). \]

\subsection{Main results}

The restricted cohomology of $\h^{\lambda}_1$  was studied in \cite{EM}. 
In this paper we compute the restricted cohomology spaces $H_*^1(\mathfrak{h}_m^{\lambda},\mathfrak{h}_m^{\lambda})$ and
$H_*^2(\mathfrak{h}_m^{\lambda},\mathfrak{h}_m^{\lambda})$ for $m>1$. With the above  notations, we can state the main results as follows:
\begin{theorem}
Over ﬁelds of characteristic $p>0$, we have
\[
\dim\ H_*^{1}(\h_m^{\lambda}(p),\h^{\lambda}_m(p))=\left\{
                             \begin{array}{ll}
                              2m^2+m+1, & \hbox{$p>2, \lambda=0$;}\\
                              2m^2+m,                                                                           & \hbox{$p>2, \lambda\neq 0$;}\\
                              2m^2-m+1,                                                                           & \hbox{$p=2, \lambda_{2m+1}=0$;}\\
                               2m^2-m, & \hbox{$p=2, \lambda_{2m+1}\neq 0$;}                            
                             \end{array}
                                                                                   \right.
\]
and
\[
\dim\ H_*^{2}(\h_m^{\lambda}(p),\h^{\lambda}_m(p))=\left\{
                             \begin{array}{ll}
                              \frac{8}{3}m^3+\frac{4}{3}m+1, & \hbox{$p>2, \lambda=0$;}\\
                              \frac{8}{3}m^3+\frac{4}{3}m,                                                                           & \hbox{$p>2, \lambda\neq 0$;}\\
                              \frac{8}{3}m^3-\frac{2}{3}m+1,                                                                           & \hbox{$p=2, \lambda_{2m+1}=0$;}\\
                               \frac{8}{3}m^3-\frac{2}{3}m, & \hbox{$p=2, \lambda_{2m+1}\neq 0$.}                            
                             \end{array}
                                                                                   \right.
\]
\end{theorem}
\begin{proof}
It follows from Propositions \ref{p1}, \ref{p2}, \ref{p3} and \ref{p4}.
\end{proof}

\begin{rem}
In fact, the results on $\h^{\lambda}_1$ in \cite{EM} satisfy the above dimensional formulas except for $H_*^{2}(\h_1^{\lambda}(p),\h^{\lambda}_1(p))$
when $p>2$.
\end{rem}

\section{Restricted cohomology $H_*^{q}(\h_m^{\lambda}(p),\h^{\lambda}_m(p))$ for $q = 1, 2$}

The results on   the ordinary cohomology spaces 
$H^1(\h_m,\h_m)$ and 
$H^{2}(\h_m,\h_m)$ are well-known and can be found in  \cite{A,M}. It is easy to observe that these results, although formulated over the complex numbers, are valid over any field, which is not surprising for nilpotent Lie algebras.

\begin{lemma}\label{H2}
For $m>1$,

(1)
$H^1(\h_m,\h_m)$ has a basis:
  \[ \left\{e^{2m+1}_{2m+1}+\sum_{i=1}^{m}e^{i}_i\right\}\bigcup\{e^{i}_{j}-e^{m+j}_{m+i}\ |\ 1\leq i,j\leq m\}\]
\[\bigcup
    \{e^{m+i}_{j}+e^{m+j}_{i},e^{i}_{m+j}+e^{j}_{m+i}\ |\ 1\le i< j\le m\}\bigcup\{e^{m+i}_i, e^{i}_{m+i}\ |\ 1\le i\le m\};\]

(2) $H^2(\h_m,\h_m)$ has a basis $A_m$ listed as follows:
  \[\{e^{ij}_{i'},e^{ij}_{j'}\ |\ 1\leq i<j\leq 2m\}\]
\[\bigcup
    \{e^{i,j'}_{j'}+e^{j,j'}_{i'},e^{ij}_{j}+e^{j,j'}_{i'}\ |\ 1\le i\neq j\le m\ \mathrm{or}\ m+1\le i\neq j\le 2m\}\]
\[\bigcup
    \{e^{kj}_{i'}+e^{ij}_{k'},e^{ki}_{j'}-e^{ij}_{k'}\ |\ 1\le k<i<j\le m\ \mathrm{or}\ m+1\le i< j<k\le 2m\}\]
\[\bigcup
    \{e^{i,k'}_{j'}-e^{ij}_{k},e^{j,k'}_{i'}+e^{ij}_{k}\ |\ k\neq i,j; 1\le i<j\le m, 1\le k\le m \ \mathrm{or}\]
\[ m+1\le i<j\le 2m, m+1\le k\le 2m\}.\]
\end{lemma}

\subsection{The Restricted Cohomology $H_*^{1}(\h^{\lambda}_m(p),\h^{\lambda}_m(p))$}

 For a restricted Lie algebra $(\g,[p])$,
 \cite[Theorem 2.1]{H}  states that
\begin{equation}\label{H1}
H^1_*(\g,\g)=\{\psi\in H^1(\g,\g)\ | \ (\ad g)^{p-1}(\psi(g))=\psi(g^{[p]})\ \mathrm{for}\  \mathrm{all}\ g\in \g\}.
\end{equation}
If $p>2$,  then  $(\ad g)^{p-1}(\psi(g))=0$ for all $g\in \h_m$ and (\ref{H1}) reduces to
\begin{align}\label{h2}
 \left\{
                             \begin{array}{ll}
H_*^{1}(\h_m^{\lambda},\h^{\lambda}_m)=
                              \{\psi\in H^1(\h_m,\h_m)\ | \ \psi(\h_m^{[p]})=0\},  & \hbox{$p>2$;}\\
H_*^{1}(\h_m^{\lambda},\h^{\lambda}_m)=  \{\psi\in H^1(\h_m,\h_m)\ | 
                                 \ [g,\psi(g)]=\psi(g^{[2]}),g\in \h_m\}
, & \hbox{$p=2$.}\\      
 \end{array}       
\right.
  \end{align}

\begin{proposition}\label{p1}
Suppose that   $p>2$ and $m>1$.

(1) If $\lambda=0$, then $H^1_*(\h^{\lambda}_m,\h^{\lambda}_m)=H^1(\h_m,\h_m)$.

(2) If $\lambda\neq 0$, then $H^1_*(\h^{\lambda}_m,\h^{\lambda}_m)=H^1(\h_m,\h_m)/\langle e^{2m+1}_{2m+1}+\sum_{i=1}^{m}e^i_i\rangle$.
\end{proposition}
\begin{proof}
From Eq.  (\ref{r}), when $p>2$,
\[\langle (\h_m^{\lambda})^{[p]}\rangle_{\F}=\left\{
                             \begin{array}{ll}
                              0, & \hbox{$\lambda=0$;}\\
                              \langle e_{2m+1}\rangle, & \hbox{$\lambda\neq 0.$}                            
                             \end{array}
                                                                                   \right.
\]
The proof is complete from Lemma  \ref{H2} (1) and Eq.  (\ref{h2}).
\end{proof}

\begin{proposition}\label{p2}
Suppose that $p=2$ and $m>1$.

(1) If $\lambda_{2m+1}=0$, then 
$H^1_*(\h^{\lambda}_m(2),\h^{\lambda}_m(2))=H^1(\h_m,\h_m)/ (\langle e^{2m+1}_{2m+1}+\sum_{i=1}^{m}e^i_i,$
\\
$e^{m+i}_{i},e^i_{m+i}\ |\ 1\le i\le m\rangle)\oplus \langle e^{2m+1}_{2m+1}+\sum_{i=1}^{m}(e^i_i+\lambda_ie^i_{m+i}+\lambda_{m+i}e^{m+i}_i)\rangle.$

(2) If $\lambda_{2m+1}\neq 0$, then 
$H^1_*(\h^{\lambda}_m(2),\h^{\lambda}_m(2))=H^1(\h_m,\h_m)/\langle e^{2m+1}_{2m+1}+\sum_{i=1}^{m}e^i_i,$
\\
$e^{m+i}_{i},e^i_{m+i}\ |\ 1\le i\le m\rangle.$
\end{proposition}
\begin{proof}
For $\psi\in H^1(\h_m,\h_m)$ and $g\in \h_m$, if we write
$g=\sum^{2m+1}_{i=1} a_i e_i$ and
\begin{align*}
\psi&=\alpha\left(e^{2m+1}_{2m+1}+\sum_{i=1}^{m}e^{i}_i\right)+\sum_{1\le i,j\le m}\alpha_{ij}\left(e^{i}_{j}-e^{m+j}_{m+i}
\right)\\
&+\sum_{1\le i< j\le m}\beta_{ij}\left(e^{m+i}_{j}+e^{m+j}_{i}\right)
+\sum_{1\le i< j\le m}\gamma_{ij}\left(e^{i}_{m+j}+e^{j}_{m+i}\right)\\
&+\sum_{i=1}^{m}\alpha_ie^{m+i}_i+\sum_{i=1}^{m}\beta_ie^{i}_{m+i},
\end{align*}
 then Eq. (\ref{r}) and Lemma \ref{H2} (1) give
                          \[\psi(g^{[2]})=   \left(\sum^{2m+1}_{i=1}a^{2}_i\lambda_i+\sum^{m}_{j=1} a_j
     a_{m+j} \right)\alpha e_{2m+1}\]
and 
\[[g,\psi(g)]=\left(-\alpha \sum_{i=1}^ma_ia_{m+i}-\sum_{i=1}^m\alpha_ia_{m+i}^2+\sum_{i=1}^m\beta_ia_i^{2}\right)e_{2m+1}.\]
From $[g,\psi(g)]=\psi(g^{[2]})$, we get
\[\alpha\lambda_{2m+1}=0,\ \alpha\lambda_i=\beta_i,\ \alpha\lambda_{m+i}=\alpha_i,\]
for $1\le i\le m$. If $\lambda_{2m+1}=0$, then $\alpha_i=\alpha\lambda_{m+i}$, $\beta_i= \alpha\lambda_i$
and
\begin{align*}
\psi&=\alpha\left(
e^{2m+1}_{2m+1}+\sum_{i=1}^{m}\left(e^{i}_i+\lambda_{m+i}e^{m+i}_i+\lambda_ie^{i}_{m+i}\right)\right)\\
&+\sum_{1\le i,j\le m}\alpha_{ij}\left(e^{i}_{j}-e^{m+j}_{m+i}
\right)+\sum_{1\le i< j\le m}\beta_{ij}\left(e^{m+i}_{j}+e^{m+j}_{i}\right)
\\
&+\sum_{1\le i< j\le m}\gamma_{ij}\left(e^{i}_{m+j}+e^{j}_{m+i}\right).
\end{align*}
If $\lambda_{2m+1}\neq 0$, then $\alpha=\alpha_i=\beta_i=0$ and
\begin{align*}
\psi&=\sum_{1\le i,j\le m}\alpha_{ij}\left(e^{i}_{j}-e^{m+j}_{m+i}
\right)+\sum_{1\le i< j\le m}\beta_{ij}\left(e^{m+i}_{j}+e^{m+j}_{i}\right)\\
&
+\sum_{1\le i< j\le m}\gamma_{ij}\left(e^{i}_{m+j}+e^{j}_{m+i}\right).
\end{align*}
\end{proof}

\subsection{The Restricted Cohomology $H_*^{2}(\h^{\lambda}_m(p),\h^{\lambda}_m(p))$}

For a restricted Lie algebra $\g$  and a restricted
$\g$-module $M$ ($g^{[p]} x = g^p x$ for all
$g\in \g$ and $x\in M$), there is a six-term exact sequence due to Hochschild \cite[p. 575]{H} that relates the ordinary and restricted 1- and 2-cohomology spaces:
\begin{diagram}[LaTeXeqno]
  \label{sixterm}
  0 &\rTo &H^1_*(\g,M)&\rTo &H^1(\g,M)&\rTo&\Hom_{\rm Fr}(\g,M^\g)&  \rTo \\
  & \rTo & H^2_*(\g,M)&\rTo &H^2(\g,M)&\rTo^H&\Hom_{\rm Fr}(\g,H^1(\g,M)). &
\end{diagram}
In the case $\g=M=\h^{\lambda}_m$,
the map
\[H: H^2(\h_m,\h_m)\to\Hom_{\rm Fr}(\h_m,H^1(\h_m,\h_m))\]
 in
(\ref{sixterm}) is given explicitly by
\[
H_\phi(g)\cdot h=
 \phi(g\wedge (\ad g)^{p-1}(h))+[g,\phi(g\wedge (\ad g)^{p-2}(h))]-\phi(g^{[p]}\wedge h),\]
where $g,h\in\h_m$ \cite{H,Viv}. 
It follows that if  
$(\phi,\omega)\in C^2_*(\h^{\lambda}_m,\h^{\lambda}_m)$, then
\begin{equation}\label{cc}
\ind^2(\phi,\omega)(h,g)=H_\phi(g)\cdot h+[h,\omega(g)].
\end{equation}
for all $g,h\in\h_m$.
Set
$$H^{[p]}_0=\{f\in\Hom_{\rm Fr}(\h_m,\F e_{2m+1})\ |\ f=\mathrm{ind}^{1}(\psi)\ \rm for\ some\ \psi\in \ker(\emph d^{1}) \}.$$
Then $H^{[p]}_0$ is the kernel of the map $\Hom_{\rm Fr}(\h_m,\F e_{2m+1})\rightarrow H^2_*(\h^{\lambda}_m,\h^{\lambda}_m)$ and the six-term exact sequence (\ref{sixterm}) decouples to the exact
sequence
\begin{diagram}[LaTeXeqno]
  \label{sixterm1}
0&\rTo &\Hom_{\rm Fr}(\h_m,\F e_{2m+1})/H^{[p]}_0 & \rTo & H^2_*(\h^{\lambda}_m,\h^{\lambda}_m)&\rTo &H^2(\h_m,\h_m)&\\
&\rTo^H&\Hom_{\rm Fr}(\h_m,H^1(\h_m,\h_m)).
\end{diagram}
Denote by $A^{[2]}_m=\{(e^{ij}_{i'},\widetilde e^{ij}_{i'}+\overline e^i_{j'}),
(e^{ij}_{j'},\widetilde e^{ij}_{j'}+\overline e^j_{i'})\ | \ 1\le i< j\le 2m\}
\bigcup
\{(\phi,\widetilde \phi)\ | \ \phi\in \left(A_m-\{e^{ij}_{i'},e^{ij}_{j'} \ | \ 1\le i<j\le 2m\}\right)\}.$

\begin{lemma}\label{rcc}
The map $H=0$ and

(1) if $p>2$, then $H^2_*(\h^{\lambda}_m,\h^{\lambda}_m)$ has a subspace isomorphic  to $H^2(\h_m,\h_m)$, spanned by
$\{(\phi,\widetilde \phi)\ | \ \phi\in A_m\}$;

(2) if $p=2$, then  $H^2_*(\h^{\lambda}_m,\h^{\lambda}_m)$ has a subspace isomorphic  to $H^2(\h_m,\h_m)$, spanned by
$A^{[2]}_m$.
\end{lemma}

\begin{proof}
Let $\phi$ be a cochain in Lemma \ref{H2} (2).   Then Eq. (\ref{r}) and  Lemma \ref{H2} (2) imply
 $\phi(g^{[p]},h)=0$ for $g,h\in\h_m$. If $p>2$, then $H_\phi=0$. From Eq. (\ref{cc}),  the cohomology classes 
$\{(\phi,\widetilde \phi)\ | \ \phi\in A_m\}$ spans a subspace of $H^2_*(\h^{\lambda}_m,\h^{\lambda}_m)$.
If $p=2$, then
$$H_{e^{ij}_{i'}}(e_k)=d^0(\delta_{ki} e_{j'}), \ H_{e^{ij}_{j'}}(e_k)=d^0(\delta_{kj} e_{i'}),\ H_\phi(e_k)=0,$$
where $\phi\in \left(A_m-\{e^{ij}_{i'},e^{ij}_{j'} \ | \ 1\le i<j\le 2m\}\right)$,
$1\le i<j\le 2m$ and $1\le k\le 2m+1$.
From Eq. (\ref{cc}),  the cohomology classes 
$A^{[2]}_m$ spans a subspace of $H^2_*(\h^{\lambda}_m,\h^{\lambda}_m)$.
\end{proof}

Since $H=0$, the exact sequence (\ref{sixterm1}) gives the splitting exact sequence
\begin{diagram}
0&\rTo &\Hom_{\rm Fr}(\h_m,\F e_{2m+1})/H^{[2]}_0 & \rTo & H^2_*(\h^{\lambda}_m,\h^{\lambda}_m)&\rTo &H^2(\h_m,\h_m)&\rTo&0, 
\end{diagram}
which follows:

\begin{corollary}\label{rcc1}
$H^2_*(\h^{\lambda}_m,\h^{\lambda}_m)\cong(\Hom_{\rm Fr}(\h_m,\F e_{2m+1})/H^{[2]}_0)\oplus H^2(\h_m,\h_m).$
\end{corollary}

\begin{proposition}
  \label{p3}
  Suppose that $p> 2$ and $m>1$.

(1) If $\lambda=0$, then
  a basis of  $H^2_*(\h_m^\lambda,\h_m^{\lambda})$ consists of the cohomology
  classes 
\[(A_m^{[p]}=\{(\phi,\widetilde \phi)\ | \ \phi\in A_m\})
    \bigcup \{(0,\overline e^i_{2m+1})\ |\ 1\le i\le 2m+1\}.\]

(2) If there is some $\lambda_i\neq 0$ for $1\le i\le 2m+1$, then a basis of $H^2_*(\h_m^{\lambda},\h_m^{\lambda})$ consists of the cohomology classes
\[(A_m^{[p]}=\{(\phi,\widetilde \phi)\ | \ \phi\in A_m\})
    \bigcup \{(0,\overline e^j_{2m+1})\ |\ 1\le j\le 2m+1\ \mathrm{and}\  j\neq i \}.\]
\end{proposition}

\begin{proof}
If $p>2$, 
$
H^{[p]}_0=\left\langle \sum_{i=1}^{2m+1}\lambda_i\overline e^{i}_{2m+1}\right\rangle
=\left\{
                     \begin{array}{ll}
0, & \hbox{$\lambda=0$;}\\
\F \overline e^i_{2m+1},& \hbox{$\lambda_i\neq 0$.}  
               \end{array}
                 \right.
$
Together with Lemma \ref{rcc} (1) and Corollary \ref{rcc1},  the proof is complete.
\end{proof}

\begin{proposition}
  \label{p4}
  Suppose that $p=2$ and $m>1$.

(1) If $\lambda_{2m+1}=0$, then
  a basis of  $H^2_*(\h_m^\lambda(2),\h_m^{\lambda}(2))$  is
$A^{[2]}_m
 \bigcup \{(0,\overline e^{2m+1}_{2m+1})\}.$

(2) If $\lambda_{2m+1}\neq 0$, then $H^2_*(\h_m^{\lambda}(2),\h_m^{\lambda}(2))\cong H^2(\h_m,\h_m)$  and a basis is $A^{[2]}_m$.
\end{proposition}

\begin{proof}
If $p=2$, 
$
H^{[2]}_0=\left\langle \lambda_{2m+1}\overline e^{2m+1}_{2m+1},\overline e^i_{2m+1}\ |\ 1\le i\le 2m\right\rangle
$\\
$
=\left\{
                     \begin{array}{ll}
\left\langle \overline e^i_{2m+1}\ |\ 1\le i\le 2m\right\rangle, & \hbox{$\lambda_{2m+1}=0$;}\\
\Hom_{\rm Fr}(\h_m,\F e_{2m+1}),& \hbox{$\lambda_{2m+1}\neq 0$.}  
               \end{array}
                 \right.
$
Together with Lemma \ref{rcc} (2) and Corollary \ref{rcc1}, the proof is complete.
\end{proof}


\begin{thebibliography}{99}
  
\bibitem{A} M. A. Alvarez.  The adjoint homology of Heisenberg Lie
  algebras.  \emph{J. Lie Theory} \textbf{25} (2015) 91--104.


%\bibitem{BaLi} W. Bai, W. D. Liu.  Cohomology of Heisenberg Lie
  %superalgebras.  \emph{J. Math. Phys.} \textbf{58} (2017), 15 pages.


% \bibitem{BoGr} S. Bouarroudj, P. Grozman, A. Lebedev, D. Leites.
%   Divided power (co)homology. Presentations of simple finite
%   dimensional modular Lie superalgebras with Cartan matrix.
%   \emph{Homology, Homotopy Appl.} \textbf{12} (2010) 237--278.


\bibitem{CaTi} L. Cagliero, P. Tirao.  The cohomology of the cotangent
  bundle of Heisenberg groups.  \emph{Adv. Math.}  \textbf{181} (2004)
  276--307.

%\bibitem{CaTi2} L. Cagliero, P. Tirao.  On the adjoint homology of
  %2-step nilpotent Lie algebras.  \emph{Bull. Austral. Math. Soc.}
  %\textbf{71} (2005) 177--182.


\bibitem{CJ} G. Cairns, S. Jambor.  The cohomology of the Heisenberg
  Lie algebras over fields of finite characteristic.
  \emph{Proc. Amer. Math. Soc.} \textbf{136} (2008) 3803--3807.

\bibitem{T} A. J. Calder\'{o}n Mart\'{i}n, C. Draper,
  C. Mart\'{i}n-Gonz\'{a}lez, J. M. S\'{a}nchez, Delgado.  Gradings
  and symmetries on Heisenberg type algebras.  \emph{ Linear Algebra
    Appl.} \textbf{458} (2014) 463--502.
    
\bibitem{ChE} C. Chevalley, S. Eilenberg, Cohomology theory of Lie
  groups and Lie algebras, \emph{Trans. Amer. Math. Soc.}  \textbf{63}
  (1948), 85-124.
    
\bibitem{EM} Q. Ehret, A. Maklouf.  Deformations and cohomology of
  restricted Lie-Rinehart algebras in positive characteristic.
  \emph{Arxive: 2305.16425}, 2023
  
%\bibitem{EFi1} T. J. Evans, A. Fialowski. Restricted central extensions
  %of restricted simple Lie algebras. \emph{Linear Algebra Appl.}
  %\textbf{513} (2017) 96-102.
    

\bibitem{EFi2} T. J. Evans, A. Fialowski.  Restricted one-dimensional
  central extensions of the restricted filiform Lie algebras
  $\mathfrak{m}^{\lambda}_{0}(p)$.  \emph{Linear Algebra Appl.}
  \textbf{565} (2019) 244--257.


  \bibitem{EFY} T. J. Evans, A. Fialowski, Yong Yang.  On the cohomology of restricted Heisenberg Lie algebras.
 \emph{Linear Algebra Appl.} \textbf{699} (2024) 295--309.
  



\bibitem{EFu}T. J. Evans, D. B. Fuchs.  A complex for the cohomology
  of restricted Lie algebras.  \emph{J. Fixed Point Theory Appl.}
  \textbf{3} (2008) 159--179.

\bibitem{Fe}J. Feldvoss. On the cohomology of restricted Lie
  algebras. \emph{Comm. Algebra} \textbf{19} (1991), 2865--2906.

%\bibitem{Fi}
%A. Fialowski. 
%An example of formal deformations of Lie algebras. 
%\emph{Deform. Theory Algebra Appl.}
%\textbf{247} (1988) 375--401.


  % \bibitem{FP} A. Fialowski, M. Penkava.  On the cohomology of Lie
  %   algebras with an invariant inner product.
  %   \emph{Algebr. Represent. Theory} \textbf{25} (2022) 1107--1131.
  
\bibitem{F} D.B. Fuchs. Cohomology of infinite dimensional Lie
  algebras. \emph{Contemporary Soviet Mathematics}, Consultants
  Bureau, New York, 1986.


\bibitem{H} G. Hochschild.  Cohomology of restricted Lie algebras.
  \emph{Amer. J. Math.} \textbf{76} (1954) 555--580.

  % \bibitem{HS} G. Hochschild, J-P. Serre.  Cohomology of Lie
  %   algebras.
  %   \emph{ Ann. Math.}  \textbf{57} (1953) 591--603.


\bibitem{J} N. Jacobson.  Lie Algebras.  \emph{John Wiley} (1962).


\bibitem{M} L. Magnin.  Cohomologie adjointe des alg\`{e}bres de Lie de
  Heisenberg.  \emph{Comm. Algebra} \textbf{21} (1993) 2101--2129.


\bibitem{S} L. J. Santharoubane.  Cohomology of Heisenberg Lie
  algebras.  \emph{Proc. Amer. Math. Soc.} \textbf{87} (1983) 23--28.

\bibitem{Se} G.B. Seligman.  Modular Lie Algebras, Ergebnisse der
  Mathematikund ihrer Grenzgebiete, Band 40, Springer 1967.


\bibitem{SK} E. Sk\"{o}ldberg.  The homology of Heisenberg Lie
  algebras over fields of characteristic two.
  \emph{Math. Proc. R. Ir. Acad.} \textbf{105A} (2005) 47--49.

%\bibitem{SF} H. Strade, R. Farnsteiner.  Modular Lie algebras and
 % their representations.  \emph{Monographs and Textbooks in Pure and
%Applied Math. Vol.116}, Marcel Dekker, Inc., New York, 1988

  % \bibitem{t} S. T\^{o}g\^{o}.  Outer derivations of Lie algebras.
  %   \emph{Trans. Amer. Math. Soc.} \textbf{128} (1967) 264--276.

 

\bibitem{Viv} F. Viviani. Restricted infinitesimal deformations of
  restricted simple {L}ie algebras.  \emph{J.  Algebra Appl.}
  \textbf{11} (5) (2012) 19 pages.






%\bibitem{Y} W. D. Liu, Y. Yang, X. K. Du.  Adjoint cohomology of
%  two-step nilpotent Lie superalgebras.  \emph{J. Lie Theory}
 % \textbf{31} (2021) 221--232.


\end{thebibliography}
\end{document}